\documentclass[a4paper,12pt]{article}
\usepackage[dvips]{epsfig}
\usepackage{amsmath,amssymb, amsbsy, amstext,amscd,amsfonts}
\usepackage{color,enumerate,euscript,graphicx,hyperref}
\usepackage{srcltx,tikz}
\input amssym.def
\input amssym.tex

\newtheorem{lem}{Lemma}[section]%
\newtheorem{theorem}[lem]{Theorem}%
\newtheorem{prop}[lem]{Proposition}%

 \def\s{\sigma} \def\t{\tau}

 \def\O{\Omega}

 \def\og{\overline G} \def\oh{\overline H}  \def\oc{\overline C}

 \def\ox{\overline X}  \def\o1{\overline 1}

 \def\ola{\overline a} 
 \def\olz{\overline z} \def\olc{\overline c} \def\olb{\overline b}
  \def\olx{\overline x} 

\def\o{\overline}   \def\olb{\overline b}
\def\di{\bigm|} \def\lg{\langle} \def\rg{\rangle}

\def\Aut{\hbox{\rm Aut\,}} \def\Inn{\hbox{\rm Inn}} \def\Syl{\hbox{\rm Syl}}
  
  \def\mod{\hbox{\rm mod }}

  \def\GL{\hbox{\rm GL}}  \def\P\GL{\hbox{\rm P\GL}}
  \def\FF{{\hbox{\sf F\kern-.43emF}}}
 \def\char{\hbox{\rm char}}

\def\Sym{\hbox{\rm Sym}}

\def\C{\hbox{\rm C}}

\def\o{\hbox{\rm o}}

\def\char{ \, {\rm char}\,}

\def\calm{\mathcal{M}}

\def\ZZ{\mathbb{Z}}  

\def\nd{\mathrel{\bigm|\kern-.7em/}} 
 \def\f{\noindent}
\def\qed{\hfill $\Box$} \def\demo{\f {\bf Proof}\hskip10pt}

\begin{document}
\begin{center}
{\bf\large  Skew Product Groups for 2-Groups of Maximal Class}
\end{center}


\begin{center}
Wenjuan Luo and Hao Yu{\small \footnotemark}\\
\medskip
 {\small
Capital Normal University,\\ School of Mathematical Sciences,\\
Beijing 100048, People's Republic of China
}
\end{center}

\footnotetext{Corresponding author: 3485676673@qq.com.
This work is supported in part by the National Natural Science Foundation of China (12071312).}

\renewcommand{\thefootnote}{\empty}
\begin{abstract}
  Skew morphisms, which generalise automorphisms for groups, provide a fundamental tool for the study of regular Cayley maps and, more generally, for finite groups with a
complementary factorisation $X = GY$,  where $Y$ is cyclic and core-free in $X$. $X$ is called the skew product group associated with $G$ and $Y$.
  In this paper, we
classify  skew product groups for the maximal class 2-groups.
  \end{abstract}

\section{Introduction}
All groups in this paper are assumed to be finite. A \textit{skew-morphism} of a group $G$ is a permutation $\s$ on $G$, having the properties that $\s(1_G)=1_G$ and there exists an integer-valued function $\pi$ on $G$ such that
$\s(gh)=\s(g)\s^{\pi(g)}(h)$ for all $g,h\in G$. $\pi$ is $\s$'s associated {\it power function}.
Note that if
$\pi(g)=1$ for all $g\in G$, then the skew-morphism $\s$ is an automorphism of $G$. Thus
skew morphisms generalise the concept of automorphisms for groups.

\smallskip
The investigation of skew-morphisms is at least  related to the following two topics.

\smallskip
(1) {\it  Group factorizations}: Use
$L_{G}:=\{L_{g}\mid g\in G\}$ to denote the left regular representation of $G$.  Then $\s$
, $L_g$ $\in \Sym(G)$. For any $g,h\in G$,  we have
\begin{equation*}
(\s L_{g})(h)=\s(gh)=\s(g)\s^{\pi(g)}(h)=(L_{\s(g)}\s^{\pi(g)})(h),
\end{equation*}
and so $\s L_{g}=L_{\s(g)}\s^{\pi(g)}$.
Therefore,  $\langle \s\rangle L_{G}\subseteq L_{G}\langle \s\rangle $. Since $|\langle \s\rangle L_{G}|=|L_{G}\langle \s\rangle|$,  we have
$\langle \s\rangle L_{G}=L_{G}\langle \s\rangle $, which implies  that $X:=L_{G}\langle\s \rangle$ is a subgroup of $\Sym(G)$, called the \textit{skew-product} of $L_{G}$ by $\s$,  see
\cite{CJT2016, ZD2016}. Moreover, one can show that $\langle \s \rangle  $  is core-free in $X$, meaning that there is no nontrivial normal subgroup of $X$ contained in $\lg \s\rg $.

Conversely, let $X$ be a finte group admitting a factorization $X=GY$
with $G\cap Y=1$ and  $Y=\langle y\rangle$ being cyclic and core-free in $X$. Then for any $g\in G$, there exists a unique $g'\in G$ and a unique $i \in \{1,2,\ldots,|Y|-1\}$ such that $yg=g'y^{i}$. This induces a permutation $\s$ on $G$ by $\s(g)=g'$, and an integer-valued function $\pi$ on $G$ by $\pi (g)=i$. Then one may check that $\s$ is a skew-morphism of $G$ with power function $\pi$.

\smallskip
(2) {\it Cayley maps}:
The concept of skew morphism was first introduced as a fundamental tool for the study of
regular {\it Cayley maps}~\cite{JS2002}. Let $G$ be a group and let $S$ be a subset of $G$ such that $1_G\not\in S, S=S^{-1}$ and $G=\langle S\rangle$. Let $\rho$ be a cycle on $S$. A Cayley map $\calm =\mathrm{CM}(G,S,\rho)$ is a 2-cell embedding
of a Cayley graph $\mathrm{Cay}(G,S)$ into an orientable closed surface such that, at each vertex $g$ of $\calm $, the local
orientation $R_g$ of the darts $(g,gx)$ incident with $g$ agrees with $\rho$ on $S$, that is, $R_g(g,gx)=(g,gx^{\rho})$ for all
$g\in G$ and $x\in S$. The automorphism group $\Aut(\calm)$ of a Cayley map $\calm$ contains a
vertex-regular subgroup induced by left multiplication of the elements of $G$ and acts semi-regularly on the darts of $\calm$. If $\Aut(\calm)$ is regular, then the map $\calm$ is called a
\textit{regular Cayley map}. It was shown by Jajcay and \v{S}ir\'{a}\v{n} that a Cayley map $\calm$ is
regular if and only if $\rho$ extends to a skew-morphism of $G$, see \cite[Theorem~1]{JS2002}. Thus the problem
of determining all regular Cayley maps of a group $G$ is equivalent to the problem of determining
all skew-morphisms of $G$ containing a generating orbit which is closed under taking inverses. Therefore, it is sufficient for us to consider skew product groups $X=GY$
with $G\cap Y=1$ and  $Y=\langle y\rangle$ being cyclic and core-free in $X$.

\smallskip
Now we are ready to recall the studying history of skew-morphisms of groups.
An interesting and important problem in this area is a determination of the skew-morphisms of a given family of groups.
The problem seems challenging because even skew-morphisms of the cyclic groups have not yet been completely determined.
For  partial results of cyclic groups, see~\cite{CJT2016,CT,DH,KN1,KN2,Kwo}.
For finite nonabelian simple group or finite nonabelian characteristically simple groups,
they were classified  in  \cite{BCV2019} and \cite{CDL}, respectively,
and for elementary abelian $p$-groups, a global structure  was characterized in \cite{DYL}.
Based on big efforts of several authors working on regular Cayley maps (see \cite{CJT2016,HKK2022,KKF2006, KMM2013,KK2017,KK2016,RSJTW2005,KK2021,WF2005,WHY2019, YWQ, YWQ2, Zhang2015,
Zhang20152,ZD2016}),
the final classification of  skew product groups of dihedral groups was given in \cite{HKK2022}.
For generalized quaternion groups, there are some partial results, see \cite{HR2022} and \cite{KO2008}.
A 2-group of order $2^n\ge 8$
 is said to be of maximal class if it has nilpotency class $n-1$.
In this paper, we shall classify skew product groups for 2-groups of maximal class.

\smallskip
Given  the \emph{skew product group} $L_G\langle \sigma\rangle$ of $L_G$ by $\sigma$,  for the purpose of this paper,  we may define  the \emph{skew product group} $X:=G\langle \sigma\rangle$ of $G$ by $\s$
as follows: every element of $X$ is uniquely written as $g\sigma^{i}$ where $g\in G$ and $i$ is a positive integer less than the order of $\sigma$; for each pair of elements $a\sigma^{i}, b\sigma^{j}\in X$, we have $(a\sigma^{i})(b\sigma^{j})=
a\sigma^{i}(b)\sigma^{\sum_{k=0}^{i-1}\pi(\sigma^{k}b)+j}$. It is straightforward to check by using the definition of the skew-morphism $\s$ that $X$ is indeed a group with operation defined above.  Sometimes, we just say $X$ a skew-product group of $G$ for short.

Throughout this paper, set $C=\lg c\rg $ and
\begin{eqnarray}\label{main0}
\begin{array}{ll}
&Q=\lg a, b\di a^{2n}=1,b^2=a^n,a^b=a^{-1}\rg\cong Q_{4n},\,n\ge 2,\\
&D=\lg a, b\di a^{n}=b^2=1, a^b=a^{-1}\rg \cong D_{2n}, \, n\ge 2.
\end{array}
\end{eqnarray}
Let $G\in\{ Q, D\}$ and let $X=X(G)=GC=\lg a, b\rg \lg c\rg$ be a group.
In Theorem~\ref{main3}, a classification of $X(Q)$ is given, provided that $C$ is core-free.
For skew product groups of $p$-groups, we have the following characterization:

\begin{theorem}\label{main1}
Let $X=GC$ be a group, where $G$ is a $p$-group and $C$ is a cyclic group such that $G\cap C=1$.
Set $C=C_1\times C_2$, where $C_1$ is the Sylow $p$-subgroup of $C$.
If $C_X=1$, then $F(X)=O_p(X)=G_1C_1$, where $G_1=O_p(X)\cap G\ne 1$ and $G_1C_1\rtimes C_2\lhd X$.
\end{theorem}

\begin{theorem}\label{main2}
Let $X=GC$ be a group, where $C$ is a cyclic group, and suppose that $G$ is a maximal class 2-group and $|G|=2^n\geq 32$. Assume that $G\cap C=1$ and that $C_X=1$. Then $X$ is a $2$-group.
\end{theorem}

\begin{theorem}\label{main3}
Let $X=GC$ be a 2-group, where $G$ is a maximal class group, $C$ is a cyclic group and $G\cap C=1$.
If $C_X=1$, then $G_X$ is $\lg a_0\rg$, $\lg a^2, b\rg$ or $G$.
\end{theorem}

\begin{theorem}\label{main4}
Let $X=GC$ be a 2-group, where $G$ is a maximal class group, $C$ is a cyclic group and $G\cap C=1$. Set $R$ is the defined relation of $G$. Then $X$ is isomorphic to one of the following groups:
\begin{enumerate}
  \item[\rm(1)] $X=\lg a,b,c|R,a^c=a^r,b^c=a^sb\rg,$
where $r^{2^m}\equiv1(2^{n-1})$, and $r^{2^{m-1}}\not\equiv1(2^{n-1})$ or
$s\frac{r^{2^{m-1}}-1}{r-1}\not\equiv0(2^{n-1})$. Moreover, if $G$ is a semidihedral 2-groups, then $2|s$;
  \item[\rm(2)]  $X=\lg a,b,c|R,a^2=a^{2r},\,c^b=a^{2s}c,\,c^a=a^{2t}b^uc^v\rg,$
where $r^{2^m}\equiv 1(\mod 2^{n-2})$, $s\sum_{l=1}^{2^m} r^{l}\equiv  0(\mod 2^{n-2})$,
either
\begin{enumerate}
  \item[\rm(2.1)] $u=0$, $r^{v-1}\equiv 1(\mod 2^{n-2})$,
  $(s+2t)r\equiv (1-r)+s\sum_{l=1}^{v}r^l(\mod 2^{n-2}),$
  $t\sum_{l=1}^{2^m}r^{l}\equiv0(\mod 2^{n-2}).
v^2\equiv1(\mod 2^m)\quad{\rm{and}}\quad
1-r\equiv tr+t\sum_{l=1}^{v}r^{l}(\mod 2^{n-2}).$
  \item[\rm(2.2)] $u=1$, $r^{v-1}+1\equiv 0(\mod 2^{n-2}).
(sr+1-r)\sum_{l=0}^{v-1}r^{l}\equiv (s+2t+1)r(\mod 2^{n-1}).
(t(1-r^{-1})+s\sum_{l=0}^{v-1}r^{l})\sum_{l=0}^{2^{m-1}-1}r^{2l}\equiv 0 (\mod 2^{n-2}).
r^2[t(1-r^{-1})+s\frac{r^v-1}{r-1}]\frac{r^{v-1}-1}{r^2-1}+2^{n-3}i\equiv 0(2^{n-2});$
\end{enumerate}
  \item[\rm(3)]  $X=\lg a,b,c|R,(a^2)^{c^2}=a^2,(c^2)^a=a^{2s}c^{-2}, (c^2)^b=a^{2u}c^2,a^c=bc^{2y}\rg,$
where $sy\equiv 1+i2^{n-3}(\mod2^{n-2})$ and $yu\equiv -1(\mod 2^{n-3})$, $i=1$ if $G$ is a generalized quaternion group and $i=0$ if $G$ is either a dihedral group or a semidihedral group.
\end{enumerate}

\end{theorem}

\section{Preliminaries}
In this section, the notation and elementary facts used in this paper are  collected.

\subsection{Notation}
In this paper, all the groups are supposed to be finite.
We set up the notation below, where $G$ and $H$ are groups, $M$ is a subgroup of $G$,
$n$ is a positive integer and $p$ is a prime number.
\begin{enumerate}
  \setlength{\itemsep}{0ex}
  \setlength{\itemindent}{-0.5em}
  \item[] $|G|$ and $\o(g)$: the order of $G$ and an element $g$ in $G$, resp.;
  \item[] $H\leq G$ and $H<G$: $H$ is a subgroup of $G$ and $H$ is a proper subgroup of $G$, resp.;
  \item[] $[G:H]$: the set of cosets of   $G$ relative to a subgroup $H$;
  \item[] $H\lhd G$ and $H\char~G$: $H$ is a normal and characteristic subgroup of $G$, resp.;
  \item[] $G'$ and $Z(G)$: the derived subgroup and the center of $G$ resp.;
  \item[] $M_G$: the core of $M$ in $G$ which is the maximal normal subgroup of $G$ contained in $M$;
  \item[] $G\rtimes H$: a semidirect product of $G$ by $H$, in which $G$ is  normal;
  \item[] $G.H$:  an extension of $G$ by $H$, where $G$ is normal;
  \item[] $C_M(G)$: centralizer of $M$ in $G$;
  \item[] $N_M(G)$: normalizer of $M$ in $G$;
  \item[] $\Syl_p(G)$: the set of all Sylow $p$-subgroups of $G$;
  \item[] $[a,b]:=a^{-1}b^{-1}ab$, the commutator of $a$ and $b$ in $G$;
  \item[] $\O_1(G)$: the subgroup $\lg g\in G\di g^p=1\rg$ of $G$ where $G$ is a $p$-group;
  \item[] $\mho_n(G)$: the subgroup $\lg g^{p^n}\di g\in G\rg$ of $G$ where $G$ is a $p$-group;
  \item[] $O_p(G)$ and $O_{p'}(G)$: the maximal $p$-subgroup and $p'$-subgroup of $G$,resp.;
  \item[] $F(X)$: the Fitting subgroup (the product of all nilpotent normal subgroups of G).
\end{enumerate}
\subsection{Elementary facts}
\begin{prop} \cite[ Theorem 11.9]{Hup} \label{maxclass}
  Let $G$ be a maximal class group and $|G|=2^n$.
  Then G is isomorphic to one of the following three groups:
\begin{enumerate}
  \item[\rm(1)]  $D_{2^n}:=\lg a,b|a^{2^{n-1}}=b^2=1,\,a^b=a^{-1}\rg,\,n\geq3$;
  \item[\rm(2)]  $Q_{2^n}:=\lg a,b|a^{2^{n-1}}=1,\,b^2=a^{2^{n-2}},\,a^b=a^{-1}\rg,\,n\geq3$;
  \item[\rm(3)]  $SD_{2^n}:=\lg a,b|a^{2^{n-1}}=b^2=1,\,a^b=a^{-1+2^{n-2}}\rg,\,n\geq4$.
\end{enumerate}
\end{prop}
\begin{lem}\label{AUT}
Let $G$ be a maximal class group and $|G|=2^n$, where $n\geq5$.
Then $\Aut(G)$ is a $2$-group.
\end{lem}

\begin{prop}\cite[Theorem 4.5]{Hup}\label{NC}
Let $H$ be the subgroup of $G$. Then $N_G(H)/C_G(H)$ is isomorphic to a subgroup of $\Aut(H)$.
\end{prop}

\begin{prop}\cite[Theorem]{Luc} \label{order}
If $G$ is a transitive permutation group of degree $n$  with a cyclic point-stabilizer,
then $|G|\le n(n-1)$.
\end{prop}

\begin{prop} \cite[Satz 1 and Satz 2]{Ito1} \label{mateabel}
Let $G=AB$ be a group, where both $A$ and $B$ are abelian subgroups of $G$. Then
\begin{enumerate}
  \item[\rm(1)]  $G$ is meta-abelian, that is, $G'$ is abelian;
  \item[\rm(2)]  if $G\ne 1$, then  $A$ or $B$ contains   a normal subgroup $N\ne 1$ of $G$.
\end{enumerate}
\end{prop}

\begin{prop} \cite[9, Kap. III, Satz 4.2(b)]{Hup}\label{fitting}
Suppose that $G$ is a solvable group and $F(G)$ is the Fitting subgroup of $G$. Then $C_G(F(G))\leq F(G)$.
\end{prop}

\begin{prop} \cite[Remark 1.2(i)]{DH}\label{cyclic}
The order of each skew morphisms of a cyclic group of order $2^n$ is equal to $2^m$ for some  $m<n$ and so the corresponding skew-product group is a bicyclic $2$-group with a core-free cyclic factor.
\end{prop}

\section{Proof of Theorem~\ref{main1}}
\begin{lem}\label{p-group}
Let $X=GC$ be a group, where $G$ is a $p$-group and $C$ is a cyclic group such that $G\cap C=1$.
Set $C=C_1\times C_2$, where $C_1\in \Syl_p(C)$.
If $C_X=1$, then $F(X)=O_p(X)=G_1C_1$, where $G_1=O_p(X)\cap G\ne 1$ and $G_1C_1\rtimes C_2\lhd X$.
\end{lem}
\demo Since $X$ is a product of two nilpotent groups, it is solvable and so  $F(X)\neq1$.
Note that $O_{p'}(X)\leq C$. Thus $O_{p'}(X)=1$ as $C_X=1$.
Then $F(X)=O_p(X)$. Let $P=GC_1\in \Syl_p(X)$.
Obviously, $O_p(X)=\bigcap_{x\in C_2}P^x$, hence $C_1\leq O_p(X)$ and so $O_p(X)=(O_p(X)\cap G)C_1$. Note that $G_1:=O_p(X)\cap G\neq 1$ as $C_X=1$. Let $\ox=X/O_p(X)=\og\oc_2$.
Observe that $O_p(\ox)=1$, which implies $F(\ox)=O_{p'}(\ox)\le \oc_2 $.
By Proposition~\ref{??}, we have $\oc_2\le C_{\ox}(F(\ox))\le F(\ox)\le \oc_2$, and therefor $O_{p'}(\ox)=\oc_2$. Thus $O_p(X)\rtimes C_2=G_1C_1\rtimes C_2\lhd X$.
\qed

\section{Proof of Theorem~\ref{main2}}
Note that both $D_4$ and $D_8$ have a skew-morphism $Z_3$.
However, when $n\ge 4$, we have the following results.
\begin{lem}
Let $X=GC$ be a group, where $C$ is a cyclic group,
and suppose that $G$ is a maximal class group and $|G|=2^n\geq 32$.
Assume that $G\cap C=1$ and that $C_X=1$. Then $X$ is a $2$-group.
\end{lem}
\demo
  The result is true for $n=5$, and so we assume that $n>5$, and we shall proceed by induction on $n$.
  Set $C=C_1\times C_2$, where $C_1\in Syl_2(C)$. Obviously, $GC_1$ is a Sylow 2-subgroup of $X$.
  By Lemma \ref{p-group}, we have $F(X)=O_2(X)$ and $O_2(X)\rtimes C_2\lhd X$.

  Assume that $O_2(X)< P$. Take $X_1:=G_1C$.
  Then $G_1:=O_2(X)\cap G< G$, and observe that $G_1$ is a cyclic group or a maximal class group.
  For the former case, $X_1/C_{X_1}$ is a $2$-group by Proposition \ref{cyclic}, and hence $C_2\char X_1\lhd X$.
  Then $C_2\lhd X$. Since $C_X=1$, we get $C_2=1$, as desired.
  For the latter case, $32\leq|G_1|=\frac{|G|}{2}$,
  by the induction hypothesis, $X_1/C_{X_1}$ is a $2$-group, and hence $C_2\char G_1C\lhd X$.
  Obviously, $C_2=1$  as $C_X=1$.

  Now assume that $O_2(X)=P$. Let $G=\lg a,b\rg$ and $C_1=\lg c_1\rg$, where $|a|=2^{n-1}$.
  Set $a_0=a^{2^{n-2}}$.
  Note that $P=GC_1=\lg a,b,c_1\rg$. Let $\Phi(P)$ be the Frattini subgroup of $P$.
  Observe that $|P:\Phi(P)|$ is either 4 or 8 as $P=\lg a,b,c_1\rg$.
  Let $\ox=X/\Phi(P)=\og \oc$.
  For the former case, $|\og|=|G/{\Phi(P)\cap G}|=2$, then $G\cap\Phi(P)=\lg a^2,b\rg$ or $\lg a\rg$.
  If $G\cap\Phi(P)=\lg a^2,b\rg$, then $\ox=(\lg \ola\rg\times \lg\olc_1\rg)\rtimes \lg \olc_2\rg$, and therefor $\lg\olc_2\rg\lhd \ox$.
  Note that $\lg a^2,b\rg C=\Phi(P)C_2\lhd X$.
  Since $32\leq|\lg a^2,b\rg|=\frac{|G|}{2}$, by the induction hypothesis, $C_2\char \lg a^2,b\rg C$, and therefor $C_2\lhd X$.
  Observe that $C_2=1$ as $C_X=1$.
  If $\Phi(P)\cap G=\lg a\rg$, then $\ox=(\lg \olb\rg\times\lg \olc_1\rg)\rtimes \lg \olc_2\rg$, and therefor $\lg\olc_2\rg\lhd \ox$. Note that $\lg a\rg C=\Phi(P)C_2\lhd X$.
  By Proposition \ref{cyclic}, $C_2\char \lg a\rg C$, and therefor $C_2\lhd X$.
  Note that $C_2=1$ as $C_X=1$.
  For the latter case, $\lg a,c_1\rg=\lg a\rg \lg c_1\rg$ is a subgroup of $X$.
  Observe that $a^2\in \Phi(P)$ and ${c_1}^2\in \Phi(P)$.
  Then $\Phi(P)=\lg a^2\rg\lg {c_1}^2\rg$.
  If $C_2\neq 1$, then $|c_1|<|a|$, and therefor $\O_{(m)}=\lg a^{2^m}\rg\neq 1$.
  Since $\lg a_0\rg\char \O_{(m)}\char P\char X$, we get $\lg a_0\rg\lhd X$.
  Since $32\leq|G/{\lg a_0\rg}|=\frac{|G|}{2}$ and $G/{\lg a_0\rg}$ is the maximal class 2-group,
  by the induction hypothesis, $\lg a_0\rg\lg c_2\rg/{\lg a_0\rg}\lhd X/{\lg a_0\rg}$,
  and therefor $\lg a_0\rg\lg c_2\rg\lhd X$.
  Note that $C_2\lhd X$ as $C_2\char\lg a_0\rg\lg c_2\rg\lhd X$, contradicting with $C_X=1$.
 \qed

\section{Proof of Theorem~\ref{main3}}
\vskip 3mm {\it Notation:}
  Recall $D_{2^n}=\lg a,b|a^{2^{n-1}}=b^2=1,\,a^b=a^{-1}\rg$,
  $Q_{2^n}=\lg a,b|a^{2^{n-1}}=1,\,b^2=a^{2^{n-2}},\,a^b=a^{-1}\rg$ and
  $SD_{2^n}=\lg a,b|a^{2^{n-1}}=b^2=1,\,a^b=a^{-1+2^{n-2}}\rg$.
  Set $n\geq 5$. Then let $X=GC$ be a 2-group where $G\in\{D_{2^n},Q_{2^n},SD_{2^n}\}$, $C=\lg c\rg\cong\ZZ_{2^m}$,
  $G\cap C=1$ and $C_X=1$. Then $X$ is a skew product group of $G$.
  By Proposition \ref{order}, we get $m<n$, that is $o(c)\leq o(a)$.
  Set $a_0:=a^{2^{n-2}}$ and $z:=c^{2^{m-1}}$.
  Recall $\Phi(X)$ is the Frattini subgroup of $X$.
  It follows from $X=\lg a,b,c\rg$, that is $d(G)\le 3$.

We prove Theorem~\ref{main3} in the following three lemma proofs.

\begin{lem}\label{G_X}
  $G_X\ne 1.$
\end{lem}
\demo
Since $X$ is 2-group, we get $Z(X)\neq 1$.
Then for any $gc^k\in Z(X)$ where $gc^k\neq 1$, we have $g\neq 1$ as $C_X=1$.
Since $[gc^k,c]=[g,c]^{c^k}=1$, we have $[g,c]=1$ and $g\in \cap_{c^i\in C}G^{c^i}=G_X$.
Thus $G_X\neq 1$.
\qed

\begin{lem}\label{ac}
If $G_X\lneqq \lg a\rg $ and  $|G_X|\ge 4$, then $\lg a\rg \lg c\rg < X$.
\end{lem}
\demo
Suppose that $\lg a\rg\lg c\rg$ is not a group. Then $X=\lg a,c\rg$ as $\lg a\rg\lg c\rg\subseteq\lg a,c\rg$ and $|\lg a\rg\lg c\rg|=\frac{X}{2}$, and thus $|\Phi(X)|=\frac{|X|}{4}$. Observe that $G<G\Phi(X)<X$ as $G<X$ and $C<X$. Then $\Phi(X)\cap G=\lg a^2,b_1\rg$ for some $b_1\in G\setminus{\lg a\rg}$ because $2\leq|G\Phi(X)/{\Phi(X)}|\leq 4$. Note that $\Phi(X)=\lg a^2,b_1\rg\lg c^2\rg$ as $c^2\in\Phi(X)$. Since $|G_X|\geq4$ and $G_X\leq\lg a\rg$, $\lg a^{2^{n-3}}\rg\char G_X\lhd X$, and so $\lg a^{2^{n-3}}\rg\lhd X$. Let $H=C_X(\lg a^{2^{n-3}}\rg)$. Note that $X/H\lesssim\Aut(\lg a^{2^{n-3}}\rg)$ and $|\Aut(\lg a^{2^{n-3}}\rg)|=2$, and therefor $|X/H|=2$. Then $\Phi(X)<H$, and so $[a^{2^{n-3}},b_1]=1$, a contradiction. Thus $\lg a\rg\lg c\rg<X$.

\qed

\begin{lem}\label{a^2}
  If $\lg a\rg\lg c\rg \le X$, then $G_X$ is either $\lg a^2,b\rg$ or $G$.
  \end{lem}
\demo
For the contrary, assume that $G_X$ is neither $G$ nor  $\lg a^2,b\rg$. Then $G_X\le \lg a\rg.$
Let $z$ be defined as above.   Pick  $a_2\in G_X$ such that $\lg c\rg\lg a_2^2\rg/ \lg a_2^2\rg $  is core-free  in $X/\lg a_2^2\rg$, but  $\lg c\rg\lg a_2\rg/ \lg a_2\rg $ has the nontrivial  core, say   $\lg c^i\rg\lg a_2\rg/ \lg a_2\rg $   in  $X/\lg a_2\rg$. Then in $\ox:=X/\lg a_2^2\rg$,  $\mho_1(\lg \ola_2\rg \times \lg \olc^i\rg )=\lg \olc^{2i}\rg\lhd \ox $, which implies $c^i=z$. In particular, $\lg a_2\rg \rtimes \lg z\rg \lhd X$.

Considering the conjufacy of $\og $ on $\lg \ola, \olz\rg \cong D_4$,   there exists an involution $\overline{a^ib}\in \og$ exchanging $\olz$ and $\ola_2\olz$ (simply we denote $\overline{a^ib}$ by $\olb$).
Since $\ox=\overline{GC}=(\lg \ola\rg\lg \olc\rg)\rtimes \lg \olb\rg$, firstly we write  $\olc^{\olb}=\ola^{s}\olc^{t}$, where $t\neq0$.
Then  $$\olc=\olc^{{\olb}^2}=(\ola^s\olc^t)^{\olb}=\ola^{-s}(\ola^s\olc^t)^t=\olc^t(\ola^s\olc^t)^{t-1},$$
that is   $(\ola^s\olc^t)^{t-1}=\olc^{1-t}$. Then we have
$$(\olc^{t-1})^{\olb}=(\olc^{\olb})^{t-1}=(\ola^s\olc^t)^{t-1}=\olc^{1-t}.$$
 If $t\neq1$, then $\olz^{\olb}\in\lg c^{t-1}\rg^{\olb}=\lg \olc^{t-1}\rg$, contradicting  to  $\olz^{\olb}=\ola_2\olz$. So $t=1$, that is  $\olc^{\olb}=\ola^s\olc$.
 Secondly, we write  $\olc^{\olb}=\olc^{t_1}\ola^{s_1}$. With the same arguments, we may get $t_1=1$ and  $\olc^{\olb}=\olc\ola^{s_1}$.
 Therefore,  we have
$$\ola^{s}\olc=\olc^{\olb}=\olc\ola^{s_1},$$
that is  $(\ola^{s})^{\olc}=\ola^{s_1}$. Clearly $\lg \ola^s\rg =\lg \ola^{s_1}\rg$, that is
$c$ normalises $\lg \ola^{s}, \olb\rg$. Then  $$\lg \ola^{s}, \olb\rg\leq  \cap_{\olc^i\in \lg \olc\rg} \og^{\olc^i}=\cap_{\olx\in \ox} \og^{\olx}={\og}_{\ox}<\lg \ola\rg ,$$ a contradiction.
\qed

\section{Classification}
To prove Theorem 1.3, set $R:=\{a^{2n}=c^m=1,\,b^2=a^n,\,a^b=a^{-1}\}.$
Then we shall deal with the five cases in Theorem 1.1  in the following five subsections, separately.
Let $A=G.\lg t\rg $ where $G\lhd A$ and $t^l=g\in G$.
Then $t$ induces an automorphism $\t$ of $G$ by conjugacy.
Recall  that  by the cyclic extension theory of groups,
this extension is valid if and only if
$$\t^l=\Inn(g)\quad  {\rm and} \quad \t(g)=g.$$

\subsection{$G\lhd X$}
\begin{lem}
Suppose that $G\lhd X$ and $C_X=1$. Then $$X=\lg a,b,c|R,a^c=a^r,b^c=a^sb\rg,$$
where $r^{2^m}\equiv0(2^{n-1})$, and $r^{2^{m-1}}\not\equiv1(2^{n-1})$ or
$s\frac{r^{2^{m-1}}-1}{r-1}\not\equiv0(2^{n-1})$. Moreover, if $G$ is a semidihedral 2-groups, then $2|s$.

\end{lem}
\demo
Since $G\lhd X$, we set $a^c=a^r$ and $b^c=a^sb$. Let $\pi\in\Aut(G)$ such that $\pi(a)=a^r$ and $\pi(b)=a^sb$. Then $o(\pi(a))=o(a)$ and $\pi^{2^m}=1$, that is $r^{2^m}\equiv1(2^{n-1})$. Note that if $G$ is a semidihedral 2-groups, then $o(b)=o(\pi(b))=o(a^sb)=2$, and so $2|s$.

Insure $\lg c\rg_X=1$: $a^z=a^{c^{2^{m-1}}}=a^{r^{2^{m-1}}}\neq a$ or $b^z=b^{c^{2^{m-1}}}=a^{s\frac{r^{2^{m-1}}}{r-1}}b\neq b$, that is $r^{2^{m-1}}\not\equiv1(2^{n-1})$ or
$s\frac{r^{2^{m-1}}-1}{r-1}\not\equiv0(2^{n-1})$.
\qed

\subsection{$G_X=\lg a^2,b\rg$}
\begin{lem}
Suppose that $G_X=\lg a^2,b\rg$. Then
$$X=\lg a,b,c|R,a^2=a^{2r},\,c^b=a^{2s}c,\,c^a=a^{2t}b^uc^v\rg,$$
where $r^{2^m}\equiv 1(\mod 2^{n-2})$, $s\sum_{l=1}^{2^m} r^{l}\equiv  0(\mod 2^{n-2})$,
either
\begin{enumerate}
  \item[\rm(1)] $u=0$, $r^{v-1}\equiv 1(\mod 2^{n-2})$,
  $(s+2t)r\equiv (1-r)+s\sum_{l=1}^{v}r^l(\mod 2^{n-2}),$
  $t\sum_{l=1}^{2^m}r^{l}\equiv0(\mod 2^{n-2}).~
v^2\equiv1(\mod 2^m)\quad{\rm{and}}\quad
1-r\equiv tr+t\sum_{l=1}^{v}r^{l}(\mod 2^{n-2}).$
  \item[\rm(2)] $u=1$, $r^{v-1}+1\equiv 0(\mod 2^{n-2}).
(sr+1-r)\sum_{l=0}^{v-1}r^{l}\equiv (s+2t+1)r(\mod 2^{n-1}).~
(t(1-r^{-1})+s\sum_{l=0}^{v-1}r^{l})\sum_{l=0}^{2^{m-1}-1}r^{2l}\equiv 0 (\mod 2^{n-2}).~
r^2[t(1-r^{-1})+s\frac{r^v-1}{r-1}]\frac{r^{v-1}-1}{r^2-1}+2^{n-3}i\equiv 0(2^{n-2}).$
\end{enumerate}
\end{lem}
\demo
$X=((\lg a^2\rg\rtimes\lg c\rg).\lg b\rg).\lg a\rg$.
Set $a^2=a^{2r},\,c^b=a^{2s}c,\,c^a=a^{2t}b^uc^v$, where $u\in\{0,1\}$.

What we should to determine the parameters $r,s,t,u$ and $v$ by analysing three extensions.

(1) $\lg a^2\rg\rtimes \lg c\rg$, where $(a^2)^c=a^{2r}$.
Set $\pi_1\in \Aut(\lg a^2\rg)$ such that $\pi_1(a^2)=a^{2r}$.
As mentioned before, this extension is valid if and only if
$\o(\pi_1(a^{2}))=\o(a^2)=2^{n-1}$ and $\pi_1^{2^m}=1$, that is
\begin{eqnarray}\label{f2.1}
r^{2^m}\equiv 1(\mod 2^{n-2}).
\end{eqnarray}

(2) $(\lg a^2\rg\rtimes\lg c\rg).\lg b\rg$, where $c^b=a^{2s}c$.
Set $\pi_2\in \Aut((\lg a^2\rg\rtimes\lg c\rg)$: $a^2\to a^{-2}$ and $c\to a^{2s}c$.
This extension is valid if and only if the following three equalities hold:

(i) $\pi_2$ preserves $(a^2)^c=a^{2r}$, as desired.

(ii) $\o(\pi_2(c))=2^m$:
$$(a^{2s}c)^{2^m}=c^{2^m}(a^{2s})^{c^{2^m}}\cdots(a^{2s})^{c}
=c^{2^m}a^{2s\sum_{l=1}^{2^m}r^l}=c^{2^m}a^{2s\sum_{l=1}^{2^m} r^{l}}=1,$$
that is
\begin{eqnarray}\label{f2.2}
s\sum_{l=1}^{2^m} r^{l}\equiv  0(\mod 2^{n-2}).
\end{eqnarray}

(iii) $\pi_2^2=\Inn(b^2):$
Since $b^2\in Z(X)$, we get $c=c^{b^2}=(a^{2s}c)^b=a^{-2s}a^{2s}c$, as desired.

(3) $((\lg a^2\rg\rtimes\lg c\rg).\lg b\rg).\lg a\rg$, where $c^a=a^{2t}b^uc^v$ and $u\in\{0,1\}$.
Set $\pi_3\in\Aut((\lg a^2\rg\rtimes\lg c\rg).\lg b\rg)$: $a^2\to a^2$, $b\to ba^2$ and $c\to a^{2t}b^uc^v$.
We divide the proof into two cases according to $u$, separately.
\vskip 3mm
{\it Case 1: $u=0$.}
\vskip 3mm
In this case, we get $c^a=a^{2t}c^v$ and $\pi_3(c)=a^{2t}c^v$.

(i)  $\pi_3$ preserves $(a^2)^c=a^{2r}$:
\begin{eqnarray}\label{f2.3}
r^{v-1}\equiv 1(\mod 2^{n-2}).
\end{eqnarray}

(ii) $\pi_3$ preserves $c^b=a^{2s}c$:
$$\begin{array}{lcl}
(c^b)^a    &=&(c^a)^{ba^2}=(a^{2t}c^v)^{ba^2}=(a^{-2t}(a^{2s}c)^v)^{a^2}\\
           &=&(a^{-2t}c^va^{2s\sum_{l=1}^{v}r^l})^{a^2}=a^{-2t}(ca^{2-2r})^va^{2s\sum_{l=1}^{v}r^l}\\
           &=&a^{-2t}c^va^{(2-2r)\sum_{l=0}^{v-1}r^l}a^{2s\sum_{l=1}^{v}r^l}\\
           &=&c^va^{-2tr^v+2(1-r^v)+2s\sum_{l=1}^{v}r^l},\\
(a^{2s}c)^a&=&a^{2(s+t)}c^v=c^va^{2(s+t)r^v},
\end{array}$$
that is
\begin{eqnarray}\label{f2.4}
(s+2t)r\equiv (1-r)+s\sum_{l=1}^{v}r^l(\mod 2^{n-2}).
\end{eqnarray}

(iii) $\o(\pi_3(c))=2^m$:
$$1=(a^{2t}c^v)^{2^m}=a^{2t\sum_{l=1}^{2^m}r^{vl}},$$
that is
\begin{eqnarray}\label{f2.5}
t\sum_{l=1}^{2^m}r^{l}\equiv0(\mod 2^{n-2}).
\end{eqnarray}

(iv) $\pi_3^2=\Inn(a^2)$:
$$ca^{2-2r}=\Inn(a^2)(c)=\pi_3^2(c)=(a^{2t}c^v)^a=a^{2t}(a^{2t}c^v)^v
=c^{v^2}a^{2tr^{v^2}+2t\sum_{l=1}^{v}r^{vl}},$$
that is
\begin{eqnarray}\label{f2.6}
v^2\equiv1(\mod 2^m)\quad{\rm{and}}\quad
1-r\equiv tr+t\sum_{l=1}^{v}r^{l}(\mod 2^{n-2}).
\end{eqnarray}

\vskip 3mm
{\it Case 2: $u=1$.}
\vskip 3mm
In this case, we get $c^a=a^{2t}bc^v$ and $\pi_3(c)=a^{2t}bc^v$.

(i)  $\pi_3$ preserves $(a^2)^c=a^{2r}$:
\begin{eqnarray}\label{f2.7}
r^{v-1}+1\equiv 0(\mod 2^{n-2}).
\end{eqnarray}

(ii) $\pi_3$ preserves $c^b=a^{2s}c$:
$$\begin{array}{lcl}
(c^b)^a    &=&(c^a)^{ba^2}=(a^{2t}bc^v)^{ba^2}=(a^{-2t}b(a^{2s}c)^v)^{a^2}\\
           &=&a^{-2t}a^{-4}b(a^{2s}ca^{2-2r})^v\\
           &=&a^{-2t-4}b(ca^{2sr+2-2r})^v\\
           &=&a^{-2t-4-2(sr+1-r)r^{-v}\sum_{l=0}^{v-1}r^{l}}bc^v\\
           &=&a^{-2t-4+2(sr+1-r)r^{-1}\sum_{l=0}^{v-1}r^{l}}bc^v,\\
(a^{2s}c)^a&=&a^{2(s+t)}bc^v,
\end{array}$$
that is
\begin{eqnarray}\label{f2.8}
(sr+1-r)\sum_{l=0}^{v-1}r^{l}\equiv (s+2t+1)r(\mod 2^{n-1}).
\end{eqnarray}

(iii) $\o(\pi_3(c))=2^m$:
$$1=(a^{2t}bc^v)^{2^m}=(a^{2t}bc^vba^{-2t}c^v)^{2^{m-1}}
=a^{2r^2(t(1-r^{-1})+s\sum_{l=0}^{v-1}r^{l})\sum_{l=0}^{2^{m-1}-1}r^{2l}},$$
that is
\begin{eqnarray}\label{f2.9}
(t(1-r^{-1})+s\sum_{l=0}^{v-1}r^{l})\sum_{l=0}^{2^{m-1}-1}r^{2l}\equiv 0 (\mod 2^{n-2}).
\end{eqnarray}

(iv) $\pi_3^2=\Inn(a^2)$:
$$ca^{2-2r}=\Inn(a^2)(c)=\pi_3^2(c)=(a^{2t}bc^v)^a=a^{2t-2}b(a^{2t}bc^v)^v=
ca^{-2r+2r^2[t(1-r^{-1})+s\frac{r^v-1}{r-1}]\frac{r^{v-1}-1}{r^2-1}}a_0^i,$$
where $i=\frac{v+1}{2}$ if $G$ is a generalized quaternion group, and $i=0$ if $G$ is not a generalized quaternion group.
Hence
\begin{eqnarray}\label{f2.10}
r^2[t(1-r^{-1})+s\frac{r^v-1}{r-1}]\frac{r^{v-1}-1}{r^2-1}+2^{n-3}i\equiv 0(2^{n-2}).
\end{eqnarray}

(4) Insure $\lg c\rg_X=1$: If $u=0$, then
$z^a=(c^{2^{m-1}})^a=(a^{2t}c^v)^{2^{m-1}}=za^{2tr\frac{r^{2^{m-1}-1}}{r-1}}~\neq z$ or $z^b=(c^{2^{m-1}})^b=(a^{2s}c)^{2^{m-1}}=za^{2sr\frac{r^{2^{m-1}-1}}{r-1}}~\neq z$, that is
$t\frac{r^{2^{m-1}-1}}{r-1}\not\equiv0(2^{n-2})$ or $s\frac{r^{2^{m-1}-1}}{r-1}\not\equiv0(2^{n-2})$. If $u=1$, then $z^a=(c^{2^{m-1}})^a=(a^{2t}bc^v)^{2^{m-1}}=za^{2r^2[t(1-r^{-1})+s\frac{r^v-1}{r-1}]
\frac{r^{2^{m-1}}-1}{r^2-1}}\neq z$ or $z^b=za^{2sr\frac{r^{2^{m-1}-1}}{r-1}}~\neq z$, that is $[t(1-r^{-1})+s\frac{r^v-1}{r-1}]\frac{r^{2^{m-1}}-1}{r^2-1}\not\equiv 0$ or $s\frac{r^{2^{m-1}-1}}{r-1}\not\equiv 0$.
\qed
\subsection{$|G_X|=2$}
\begin{lem}\label{G_X=2}
Suppose that $n\ge 5$ and $|G_X|=2$.
Then $X=KC=((\lg a^2,b\rg\rtimes \lg c^2\rg).\lg a\rg).\lg c\rg$, where $K=G\lg c^2\rg$. Moreover, $K'=\lg a^2\rg\times \lg c_1^2\rg$, $G_K=\lg a^2,b\rg$ and $\lg a^2\rg\rtimes\lg c_1\rg\lhd X$.
\end{lem}

\demo
Suppose that $G_X=\lg a_0\rg$.
Consider the faithful permutation representation of $X$ on $[X:G]$. Since if $M,~M^c\leq G$, then $|M|\le 4$, and so $Gc^{-1}G$ contains at least $\frac{|G|}4=2^{n-2}$ coset of $G$, and $1+2^{n-2}\le |[X:G]|=2^m$.
Note that $m=n-1$ as $n-2<m<n$. Then $o(a)=o(c)$.

Since $X/G_X=X/\lg a_0\rg=(G/\lg a_0\rg) (C\lg a_0\rg/\lg a_0\rg)$ and $|G/\lg a_0\rg|=|C\lg a_0\rg/\lg a_0\rg|$,  the core of $C \lg a_0\rg/\lg a_0\rg$ in $X/\lg a_0\rg $ is $\lg z\rg\lg a_0 \rg$.  Set $\ox=X/\lg a_0,z\rg=\og\oc$, where  $\og\cong D_{2^{n-1}}$ and $\oc$ is core-free. Let $\oh=\og_{\ox}$, with $\lg a_0,z\rg\leq H$, and note that $H\lhd X$. Observe that $\og_{\ox}$ is either dihedral or cyclic. Then we have the following two cases:

\vskip 3mm
{\it Case 1: $\og_{\ox}$ is dihedral.}
\vskip 3mm
 In this case, $\og_{\ox}$ is either $\og$ or $\lg \ola^2,\olb\rg$. Then $H=\lg a,b\rg\rtimes\lg z\rg$ or $H=\lg a^2,b\rg\rtimes \lg z\rg$.
Note that $\lg a^4\rg\leq H'< G$ and $H'\char H\lhd X$, and thus $\lg a^4\rg\leq H'\leq G_X\lhd X$. Since $G_X=\lg a_0\rg$, we have $a^4=a_0$, and so $n=4$, contradict to $n\geq 5$.

\vskip 3mm
{\it Case 2: $\og_{\ox}$ is cyclic.}
\vskip 3mm
Assume that $\og_{\ox}=\lg \ola^i\rg$ is cyclic. By Lemma \ref{ac}, we get $\og_{\ox}=\lg\ola_1\rg$, and therefor $H:=\lg a_1\rg\rtimes \lg z\rg\lhd X$.
   Observe that $H\cong D_8$ or $H\cong \ZZ_4\times \ZZ_2$.
 If $H\cong D_8$, then $\lg a_1\rg\char H\lhd X$, which implies $\lg a_1\rg\leq G_X$, a contradiction.
 Then $H=\ZZ_4\times \ZZ_2$. Note that $X/{C_X(\lg a_0\rg\times\lg z\rg)}\cong \ZZ_2$ and $c\in C_X(\lg a_0\rg\times\lg z\rg)$. By Lemma \ref{a^2}, $\lg a\rg\lg c\rg$ is not a group, and hence $a\not\in C_X(\lg a_0\rg\times\lg z\rg)$.

Set $a_1^c=a_1^iz$ and $z^a=a_0z$, where $i\in\{1,-1\}$. Note that $\lg a^2,c^2\rg\leq C_X(H)$ as $a_1^{c^2}=a_1^{i^2}=a_1$ and $z^{a^2}=z$. Since $[a_1,b]\neq 1$ and $[a,z]\neq 1$, we have $G/{C_X(H)\cap G}\cong\ZZ_2\times\ZZ_2$. If $X=GC_X(H)$, then $\lg a_1\rg\lhd X$, a contradiction. Thus $GC_X(H)<X$ and $\frac{|X|}{|C_X(H)|}\geq 2^3$. Note that $X/{C_X(H)}\lessapprox \Aut(H)\cong D_8$. Hence           $C_X(H)=\lg a^2\rg\lg c^2\rg\lhd X$ and $X/C_X(H)\cong D_8$. Let $K:=GC_X(H)=G\lg c^2\rg$. Note that $a_1\in G_K$, by Theorem \ref{main3}, $G_K=G$ or $G_K=\lg a^2,b\rg$ for some $b\in G\setminus{\lg a\rg}$. If $G_K=G$, then $\lg a^2\rg \leq K'\leq G$, and so $K'\lhd X$, and hence $\lg a^2\rg\leq G_X$, a contradiction. Thus $G_K=\lg a^2,b\rg$ and $X=((\lg a^2,b\rg\rtimes\lg c^2\rg).\lg a\rg)\lg c\rg$. Note that $\lg a^2\rg\lhd K$ as $\lg a^2\rg\char G_K\lhd K$. Obviously, $\lg a^2\rg\lhd C_X(H)$, and hence $C_X(H)'\leq\lg a^2\rg$. Note that $C_X(H)'\lhd X$ as $C_X(H)'\char \C_X(H)\lhd X$. Then $C_X(H)'\leq\lg a_0\rg$. Thus $[a^2,c^4]=1$. Since $K/\lg a^2\rg\rtimes\lg c^2\rg\cong \ZZ_2^2$ and $G<K$, this means that $K'\leq\O(K)<\lg a^2\rg\lg c^2\rg$, and so we set $K'=\lg a^2\rg\times \lg c^{4j}\rg$ for some integer $j$. Since $G_X=\lg a_0\rg$ and $\O_{o(c^{4j})}(K')\char K\lhd X$, therefor $\O_{o(c^{4j})}(K')=\lg a_0\rg$, and so $K'=\lg a^2\rg\times \lg c^4\rg$.
\qed

\begin{lem}\label{GX2}
Suppose that $|G_X|=2$. Then
$$X=\lg a,b,c|R,(a^2)^{c^2}=a^2,(c^2)^a=a^{2s}c^{-2}, (c^2)^b=a^{2u}c^2,a^c=bc^{2y}\rg,$$
where $sy\equiv 1+i2^{n-3}(\mod2^{n-2})$ and $yu\equiv -1(\mod 2^{n-3})$, $i=1$ if $G$ is a generalized quaternion group and $i=0$ if $G$ is either a dihedral group or a semidihedral group.
\end{lem}
\demo
Since $X/\lg a^2\rg\rtimes \lg c^2\rg=\lg\ola,\olb\rg\lg\olc\rg\cong D_8$, we can choose $\olb$
such that  the form of $X/M_X$ is the following:
$\ola^{\olc}=\olb$ and $\olb^{\olc}=\ola.$
Set $c_1:=c^2$ and $a_1:=a^2$.
Noting $\lg a^2,b\rg\lhd \lg a,b,c_1\rg$, we can set
$$a_1^{c_1}=a_1^r,c_1^a=a_1^sc_1^t, c_1^b=a_1^{u}c_1\quad{\rm{and}}\quad
a^c=bc_1^y.$$
Then one can check $b^c=a^{1-2sr}c_1^{1-t-y}$.
Set $H:=\lg a_1\rg\rtimes\lg c_1\rg$. Then $H'\leq\lg a_1\rg$.
Since $H'\char H\lhd X$ and $|G_X|=2$, we get $H'\leq\lg a_0\rg$, which implies $c_1^2\in C_X(a_1)$.
If $y\equiv0(\mod 2^{n-1})$, then $\o(a)=\o(a^c)=\o(b)$ is either 2 or 4, which implies $|G|=4$ or 8, a contradiction. Therefore, $y\not\equiv 0(\mod 2^{n-1})$.

What we should to determine the parameters $r,s,t,u$ and $v$ by analysing three extensions.

(1) $\lg a_1\rg\rtimes \lg c_1\rg$, where $a_1^{c_1}=a_1^r$.
Set $\pi_1\in \Aut(\lg a^2\rg)$ such that $\pi_1(a^2)=a^{2r}$.
As mentioned before, this extension is valid if and only if
$\o(\pi_1(a^{2}))=\o(a^2)=2^{n-1}$ and $\pi_1^{2}=1$, that is
$r$ is either 1 or $1+2^{n-3}$.

(2) $(\lg a_1\rg\rtimes\lg c_1\rg).\lg b\rg$, where $c_1^b=a_1^{u}c_1$.
Set $\pi_2\in \Aut((\lg a_1\rg\rtimes\lg c_1\rg)$: $a_1\to a_1^{-1}$ and $c_1\to a_1^{u}c_1$.
Note that $b^2\in\lg a_0\rg\leq Z(X)$,
thus one can check that $\pi_2$ preserves $a_1^c=a_1^{r}$, $\o(\pi_2(c))=2^m$ and $\pi_2^2=\Inn(b^2)$.

(3) $((\lg a_1\rg\rtimes\lg c_1\rg).\lg b\rg).\lg a\rg$, where $c_1^a=a_1^sc_1^t$.
Set $\pi_3\in\Aut(\lg a_1\rg\rtimes\lg c_1\rg).\lg b\rg)$: $a_1\to a_1$, $c_1\to a_1^sc_1^t$ and
$b\to ba_1$ or $ba_1a_0$.

(i)  $\pi_3$ preserves $(a^2)^c=a^{2r}$, as desired.

(ii) $\pi_3$ preserves $c_1^b=a_1^uc_1$, that is $(a_1^sc_1^t)^{ba}=a_1^uc_1$:
$$\begin{array}{lcl}
a_1^uc_1&=&(a_1^sc_1^t)^{ba}\\
        &=&a_1^{-s}c_1^{t^2}a_1^{(u+s)\sum_{l=1}^{t}r^{l}},
\end{array}$$
that is
\begin{eqnarray}\label{f3.1}
t^2\equiv1(\mod 2^{n-2})\quad{\rm{and}}\quad
(u+s)(r+1)\frac{t-1}{2}\equiv 0(\mod 2^{n-2}).
\end{eqnarray}

(iii) $\o(\pi_3(c_1))=2^{n-2}$:
$$1=(a_1^sc_1^t)^{2^{n-2}}=a_1^{s\sum_{l=1}^{2^{n-2}}r^{l}},$$
that is
\begin{eqnarray}\label{f3.2}
s\sum_{l=1}^{2^{n-2}}r^{l}\equiv0(\mod 2^{n-2}).
\end{eqnarray}

(iv) $\pi_3^2=\Inn(a_1)$: Recall that $\Inn(a_1)(c_1)=c_1a_1^{1-r}$.
$$c_1a_1^{1-r}=\Inn(a^2)(c_1)=\pi_3^2(c_1)=(a_1^sc_1^t)^a
=c^{t^2}a^{sr+s\sum_{l=1}^{t}r^{l}},$$
that is
\begin{eqnarray}\label{f3.3}
t^2\equiv1(\mod 2^{n-2})\quad{\rm{and}}\quad
1-r\equiv sr+s\sum_{l=1}^{t}r^{l}(\mod 2^{n-2}).
\end{eqnarray}

(4) $((\lg a^2,b\rg\rtimes \lg c^2\rg).\lg a\rg).\lg c\rg$, where $a^c=bc_1^y$
and $b^c=a^{1-2sr}c_1^{1-t-y}$.
Set $\pi_4\in \Aut(\lg a,b,c_1\rg) :$
$a\to bc_1^y$, $b\to a^{1-2sr}c_1^{1-t-y}$ and $c_1\to c_1$.
Let $i=1$ if $G$ is a generalized quaternion group and
$i=0$ if $G$ is either a dihedral group or a semidihedral group.
We need to carry out the following  seven steps:

(i)  $\o(\pi(a))=2^{n-1}$: Since $a_0\in Z(X)$, we only show $(a^c)^{2^{n-2}}=a_0$
$$(bc^{2w})^{2^{n-2}}=a_1^{u2^{n-3}\sum_{l=1}^{y}r^l}=a_1^{2^{n-3}},$$
that is
\begin{eqnarray}\label{f3.4}
u\sum_{l=1}^{y}r^l\equiv1(\mod2),
\end{eqnarray}
which implies that both $u$ and $y$ are odd.

(ii) $\pi_4$ preserves $a_1^{c_1}=a_1^r$:

$$(a_1^{c_1})^c=c^{2y}a^{u\sum_{l=1}^{y}r^l+i2^{n-3}}\quad{\rm{and}}\quad
(a_1^{r})^c=c^{2yr}a^{u\sum_{l=1}^{y}r^l+i2^{n-3}},$$
that is
\begin{eqnarray}\label{f3.5}
2y(r-1)\equiv 0(\mod2^{n-2}).
\end{eqnarray}

(iii) $\pi_4$ preserves $c_1^a=a_1^{s}c_1^{t}$:
$$(c_1^a)^c=a_1^{ur}c_1\quad{\rm{and}}\quad
 (a_1^{s}c_1^{t})^c=a_1^{s(ur\sum_{l=1}^{y}r^l+i2^{n-3})}c_1^{t+2ys},$$
that is
\begin{eqnarray}\label{f3.6}
ur\equiv s(ur\sum_{l=1}^{y}r^l+i2^{n-3})(\mod2^{n-2})\quad{\rm{and}}\quad
1\equiv t+2ys(\mod 2^{n-2}).
\end{eqnarray}

(iv) $\pi$ preserves $c_1^b=a_1^{u}c_1$:
$$(c_1^b)^c=c_1^{t}a_1^{sr+2(sr-1)(1-r)}\quad{\rm{and}}\quad
 (a_1^{u}c_1)^c=c_1^{2yu}a_1^{u(ur\sum_{l=1}^{y}r^l+i2^{n-3})}c_1,$$
that is,
\begin{eqnarray}\label{f3.7}
t\equiv 1+2yu(\mod2^{n-2})\quad{\rm{and}}\quad
sr\equiv u^2\sum_{l=1}^{y}r^l+i2^{n-3}(\mod2^{n-2}).
\end{eqnarray}

(v) $\pi^2=\Inn(c_1)$:
Recall $\Inn(c_1)(a)=a^{1-2sr}c_1^{1-t},\,\Inn(c_1)(a_1)=a_1^{r}$
and $\Inn(c_1)(b)=a_1^{ur}b$.
$$a^{1-2sr}c^{2-2t}=\Inn(c_1)(a)=\pi^2(a)=b^cc^{2y}=a^{1-2sr}c^{2-2t-2y+2y},$$
as desired;
$$a_1^{r}=\Inn(c_1)(a^2)=\pi^2(a^2)=(c_1^{2y}a_1^{ur\sum_{l=1}^{y}r^l+i2^{n-3}})^c
 =c_1^{2y+2yur\sum_{l=1}^{y}r^l}a_1^{(u\sum_{l=1}^{y}r^l)^2},$$
that is
\begin{eqnarray}\label{f3.8}
2y+2yur\sum_{l=1}^{y}r^l\equiv 0(\mod 2^{n-2})\quad{\rm{and}}\quad
r\equiv (u\sum_{l=1}^{y}r^l)^2 (\mod 2^{n-2});
\end{eqnarray}
and
$$a^{2ur}b=\Inn(c_1)(b)=\pi^2(b)=(a^{1-2sr}c_1^{1-t-y})^c=
a_1^{-su\sum_{l=1}^{y}r^l+i2^{n-3}}c_1^{-2sy}bc_1^{1-t},$$
as desired.

Now we are  ready to determine the parameters by summarizing  Eq(\ref{f3.1})-Eq(\ref{f3.8}).

By Eq(\ref{f3.6}) and Eq(\ref{f3.7}), we get $2(u+s)\equiv0(\mod2^{n-2})$ and $s^2r\equiv u^2(\mod2^{n-2})$.
Noting that $s^2\equiv u^2(\mod2^{n-2})$ and $u$ is odd, we get that $s$ is odd and $r=1$.
Inserting $r=1$ into Eq(\ref{f3.1})-Eq(\ref{f3.8}), we get $t=-1$ by Eq(\ref{f3.3}).
Then we get $sy\equiv 1+i2^{n-3}(\mod2^{n-2})$ by Eq(\ref{f3.6}) and
$yu\equiv -1(\mod 2^{n-3})$ by Eq(\ref{f3.7}).

\qed


{\small \begin{thebibliography}{99}

\bibitem{AK2009}  B. Amberg,  and L. Kazarin,  Factorizations of groups and related topics. {\it Sci. China Ser. A} {\bf 52} (2009)(2),  217--230.

\bibitem{BCV2019} M. Bachrat\'y, M. Conder and G. Verret,  Skew-product groups for monolithic groups, arXiv:1905.00520v1, 2019.

\bibitem{CDL} J.Y. Chen, S.F. Du and C.H. Li,  Skew-morphisms of nonabelian characteristically simple groups, {\it J. Combin. Theory Ser. A} {\bf 185} (2022), paper No. 105539, 17 pp.

\bibitem{CJT2016} M. Conder, R. Jajcay and T. Tucker,
 Cyclic complements and skew-morphisms of groups, {\it J. Algebra} {\bf 453} (2016), 68--100.

\bibitem{CT} M. Conder and R. Tucker, Regular Cayley maps for cyclic groups,
{\it Trans. Amer. Math. Soc.} {\bf366} (2014), 3585--3609.

\bibitem{D1961}J. Douglas, On the supersolvability of bicyclic groups,
{\it Proc. Nat. Acad. Sci. U.S.A.} {\bf 47} (1961), 1493--1495.

\bibitem{DMM2008} S.F. Du, A. Malni\v{c} and D. Maru\v{s}i\v{c},
Classification of 2-arc-transitive dihedrants,
{\it Journal of Combinatorial Theory} {\bf 98.6} (2008), 1349--1372.

\bibitem{DH} S.F. Du and K. Hu, Skew-morphisms of cyclic $2$-groups,
{\it J. Group Theory}  {\bf 22} (2019)(4), 617--635.

\bibitem{DYL} S.F. Du, W.J. Luo, H. Yu and J.Y. Zhang,
Skew-morphisms of elementary abelian $p$-groups, arXiv:2205.07734, 2022.

\bibitem{G1952} W. Gasch\"{u}tz, Zur Erweiterunstheorie endlicher Gruppen,
{\it J. Math.} {\bf 190} (1952), 93--107.

\bibitem{H1959}M. Hall, The Theory of Groups, {\it Macmillan} 1959.

\bibitem{HKK2022} K. Hu, I. Kov\'acs and Y. S. Kwon,
Classification of skew morphisms of dihedral groups, {\it J. Group Theory} (2022), https://doi.org/10.1515/jgth-2022-0085.

\bibitem{HR2022} K. Hu and D.Y. Ruan, Smooth skew morphisms of dicyclic groups,
{\it J. Algebraic Combin} {\bf 56} (2022), 1119--1134.

\bibitem{Hup}B. Huppert, {\it Endliche Gruppen} I, Springer, Berlin,(1967).

\bibitem{Ito1} N.~It\^o, \"{U}ber das Produkt von zwei abelschen Gruppen,
{\sl Math.~Z.} {\bf 62} (1955), 400--401.

\bibitem{JS2002}R. Jajcay and J. \v{S}ir\'{a}\v{n}, Skew-morphisms of regular Cayley maps,
{\it Disc. Math.} {\bf 224}(2002), 167--179.

\bibitem{Ke1961} O.H. Kegel, Produkte nilpotenter Gruppen,
{\it Arch Math} {\bf 12} (1961), 90--93.

\bibitem{KN1}I. Kov\'acs and R. Nedela, Decomposition of skew-morphisms of cyclic groups,
{\it Ars Math. Contemp.} {\bf 4} (2011), 329--349.

\bibitem{KN2}I. Kov\'acs and R. Nedela, Skew-morphisms of cyclic $p$-groups,
{\it J. Group Theory}  {\bf 20} (2017)(6), 135--154.

\bibitem{KMM2013}I. Kov\'acs, D. Maru\v{s}i\v{c} and M.E. Muzychuk, On $G$-arc-regular
dihedrants and regular dihedral maps, {\it J. Algebraic Combin.} {\bf 38} (2013), 437--455.

\bibitem{KK2017} I. Kov\'acs and Y. S. Kwon,
Classification of reflexible Cayley maps for dihedral groups,
{\it J. Combin. Theory Ser. B} {\bf127} (2017), 187--204.

\bibitem{KK2016} I. Kov\'acs and Y.S. Kwon, Regular Cayley maps on dihedral groups
with smallest kernel, {\it J. Algebraic Combin.} {\bf 44} (2016), 831--847.

\bibitem{KK2021} I. Kov\'acs and Y.S. Kwon, Regular Cayley maps for dihedral groups,
{\it J. Combin. Theory Ser. B} {\bf 148} (2021), 84--124.

\bibitem{Kwo}Y.S. Kwon, A classification of regular $t$-balanced
Cayley maps for cyclic groups, {\it Disc. Math.} {\bf 313} (2013), 656--664.

\bibitem{KKF2006} J.H. Kwak, Y.S. Kwon and R. Feng,
A classification of regular t-balanced Cayley maps on dihedral groups,
{\it European J. Combin.} {\bf 27} (2006), 382--393.

\bibitem{KO2008} J.H. Kwak and J.M. Oh,
A classification of regular t-balanced Cayley maps on dicyclic groups,
{\it European J. Combin.} {\bf 29} (2008), 1151--1159.

\bibitem{LX2022} C.H. Li and B.Z. Xia,
Factorizations of almost simple groups with a solvable factor, and Cayley graphs of solvable groups.
{\it Mem. Amer. Math. Soc.} {\bf 279} (2022)(1375), v+99 pp.

\bibitem{LPS1990} M.W. Licheck, C. E. Prager and J. Saxl,
The maximal factorizationof the finite simple groups and their automorphism groups,
{\it Mem. Amer: Math.Soc.} {\bf 432} (1990), 1--151.

\bibitem{Luc} A. Lucchini, On the order of transitive permutation groups with cyclic point-stabilizer,
{\it Atti. Accad. Naz. Lincei CI. Sci. Fis. Mat. Natur. Rend. Lincei (9) Mat. Appl.} {\bf 9} (1998), 241--243.

\bibitem{M1974} V.S. Monakhov,
The product of two groups, one of which contains a cyclic subgroup of index $\leq$ 2,
{\it Mathematical Notes of the Academy of Sciences of the Ussr}  {\bf 16.2} (1974), 757--762.

\bibitem{RSJTW2005} B. Richter, J. \v{S}ir\'a\v{n}, R. Jajcay, T. Tucker, and M. Watkins.
Cayley maps, {\it J. Combin. Theory Ser. B} {\bf 95} (2005), 189--245.

\bibitem{S1982} M. Suzuki, Group Theory. I, {\it Springer} 1982.

\bibitem{WF2005} Y. Wang and R. Feng,
Regular Cayley maps for cyclic, dihedral and generalized quaternion groups,
{\it Acta Math. Sin. (Engl. Ser.)} {\bf 21} (2005), 773--778.

\bibitem{WHY2019} N.E. Wang, K. Hu, K. Yuan, and J.Y. Zhang,
Smooth skew morphisms of dihedral groups,
{\it Ars Math. Contemp.} {\bf 16} (2019), 527--547.


\bibitem{W1964} H. Wielandt, Finite permutation groups. {\it Academic Pr.} 1964.

\bibitem{Wie1958} H. Wielandt, \"{U}ber Produkte von nilpotenter Gruppen,
 {\it Illinois J. Math} {\bf 2} (1958), 611--618.


\bibitem{W1999}J.S. Wilson, Products of Groups, {\it Oxford Mathematical Monographs} 1999.

\bibitem{YWQ} K. Yuan, Y. Wang and H.P. Qu, Classification of regular balanced Cayley maps of minimal non-abelian metacyclic groups, {\it Ars Math. Contemp.} {\bf 14} (2018) 433--443.

\bibitem{YWQ2} K. Yuan, Y. Wang and H.P. Qu, Regular balanced Cayley maps on nonabelian metacyclic groups of odd order, {\it Art Discrete Appl. Math.} {\bf 3} (2020) P1.05.

\bibitem{Zhang2015} J.Y. Zhang, A classification of regular Cayley maps with trivial
Cayley-core for dihedral groups, {\it Discrete Math.} {\bf 338} (2015), 1216--1225.

\bibitem{Zhang20152} J.Y. Zhang, Regular Cayley maps of skew-type $3$ for dihedral groups,
{\it Discrete Math} {\bf 388} (2015), 1163--1172.

\bibitem{ZD2016} J.Y. Zhang and S.F. Du,  On the skew-morphisms of dihedral groups,
{\it J. Group Theory} {\bf 19} (2016), 993--1016.

\end{thebibliography} }
\end{document}